\documentclass{article}

\usepackage{amsmath, amssymb, latexsym}

\newtheorem{thm}{Theorem}[section]
\newtheorem{lma}{Lemma}[section]
\newtheorem{rmk}{Remark}[section]
\newtheorem{prop}{Proposition}[section]
\newtheorem{corol}{Corollary}[section]
\newtheorem{deff}{Definition}[section]

\def\H{\mathbb{H}}
\def\R{\mathbb{R}}\def\Lap{\triangle}
\def\bX{\bar{X}}
\def\bg{\bar{g}}
\def\hg{\hat{g}}
\def\pf{\noindent{\em Proof: }}
\def\stop{\hfill$\Box$}

\begin{document}
\title{Ricci Curvature Rigidity for Weakly Asymptotically Hyperbolic Manifolds}
\date{October, 2003}
\author{
Vincent Bonini\footnote{Department of Mathematics, UCSC, Santa
Cruz, CA 95064, vbonini@math.ucsc.edu}, Pengzi
Miao\footnote{Mathematical Science Research Institute, Berkeley,
CA 94720, pengzim@msri.org} and Jie Qing \footnote{Department of
Mathematics, UCSC, Santa Cruz, CA 95064, qing@math.ucsc.edu}}
\maketitle
\begin{abstract}
A rigidity result for weakly asymptotically hyperbolic manifolds
with lower bounds on Ricci curvature is proved without assuming
that the manifolds are spin. The argument makes use of a
quasi-local mass characterization of Euclidean balls from
\cite{Miao} \cite{S_T} and eigenfunction compactification ideas
from \cite{Qing}.
\end{abstract}

\section{Introduction}

Rigidity questions for asymptotically hyperbolic manifolds have
been studied by many authors under various assumptions. In
\cite{Min-Oo}, Min-Oo proved a scalar curvature rigidity theorem
for manifolds which are spin and are asymptotic to the hyperbolic
space in a strong sense. In \cite{A_D}, Andersson and Dahl
improved the scalar curvature rigidity for asymptotically locally
hyperbolic spin manifolds. They also established the rigidity for
conformally compact Einstein manifolds with spin structure. More
recent related works are in \cite {C_H}, \cite{Wang} and
\cite{Zhang}. It is interesting to ask whether the spin structure
is necessary to assure the rigidity. In \cite{listing}, Listing
was able to obtain a non-spin rigidity at the expense of replacing
scalar curvature bound by sectional curvature bound. Very
recently, in \cite{Qing}, Qing established the rigidity for
conformally compact Einstein manifolds of dimension less than $7$
without assuming spin structure. The proof in \cite {Qing} uses
conformal compactifications by positive eigenfunctions and the
classic positive mass theorem proved by Schoen and Yau
\cite{Schoen_Yau} for asymptotically flat manifolds. Based on
ideas in \cite{Qing} combined with a quasi-local mass
characterization of Euclidean balls in \cite{Miao}, in this paper
we prove a Ricci curvature rigidity theorem for weakly
asymptotically hyperbolic manifolds.

\begin{thm} \label{Rigidity}
Let $(X^{n+1}, g)$ be a weakly asymptotically hyperbolic manifold
of order $C^{3, \alpha}$. Assume that $(X, g)$ has the standard
round sphere $(S^n, [h_0])$ as its conformal infinity and
satisfies $Ric(g) \geq -n g$. Let $r$ be the special defining
function such that
\begin{equation} \label{gexpan}
g = \frac{1}{\sinh^2(r)}\{dr^2 + g_r\}
\end{equation}
in a neighborhood of $\partial X$ and $g_0 = h_0$ . Then, if $2 \leq n \leq 6$ and
\begin{equation} \label{decayassum}
Tr_{g_r}(\frac{d}{dr} g_r) \in \Lambda^s_{0, \beta}(X)
\end{equation}
for some $s > 1$, $(X, g)$ is isometric to the hyperbolic space $\H^{n+1}$.
\end{thm}

We remark that Theorem \ref{Rigidity} may be compared with the
corresponding result in the asymptotically flat case established
by Bartnik in \cite{Bartnik}. The paper is organized as follows.
In Section $2$, we introduce notations and definitions. In Section
$3$, we recall some analytic and geometric preliminaries. In
Section $4$, we perform the conformal compactification and prove
Theorem \ref{Rigidity}. We conclude the paper by comparing our
result to scalar curvature rigidity for asymptotically hyperbolic
manifolds in \cite{A_D}, \cite{C_H}, \cite{Min-Oo}, \cite{Wang}
and \cite{Zhang}.

\section{Weakly Asymptotically Hyperbolic Manifolds}

In this section we define our terms and introduce the function
spaces that we will be working with. Throughout this paper, smooth
will always mean $C^{\infty}$.

A smooth Riemannian metric $g$ in the interior $X^{n+1}$ of a
smooth compact manifold $\bX$ with boundary is said to be {\em
conformally compact of order} $C^{m, \alpha}$ if $\bar{g}= \rho^2
g$ extends as a $C^{m, \alpha}$ metric on $\bX$, where $\rho$ is a
smooth defining function for $M^n=\partial X$ in $\bX$ in the
sense that $\rho > 0$ in $X$ and $\rho = 0$, $d \rho \neq 0$ on
$M$. The metric $\bg$ restricted to $TM$ induces a metric $\hg$ on
$M$ which rescales upon change in defining function. Therefore a
conformally compact $(X^{n+1}, g)$ defines a conformal structure
on $M$. We call $(M, [\hg])$ the {\em conformal infinity of} $(X,
g)$. When $m + \alpha \geq 2$, a straightforward computation as in
\cite{Mazzeo} shows that the sectional curvatures of $g$
approaches $-|d \rho|^2_{\bg}$ at $M$. Accordingly, we have the
following definition for weakly asymptotically hyperbolic
manifolds.

\begin{deff} For a complete manifold $(X^{n+1}, g)$, we say the metric
$g$ is {\em weakly asymptotically hyperbolic of order} $C^{m,
\alpha}$ if $g$ is conformally compact of order $C^{m, \alpha}$,
$m + \alpha \geq 2$ and $|d \rho|^2_{\bg}=1$ along $M$.
\end{deff}

To illustrate the difference between weakly asymptotically
hyperbolic and asymptotically hyperbolic we recall, for instance,
the following definition from \cite{Wang}.

\begin{deff} A weakly asymptotically hyperbolic manifold $(X^{n+1}, g)$
is called asymptotically hyperbolic if it satisfies:

(1) The conformal infinity is the round sphere one $(S^n,
[h_0])$.

(2) For a geodesic defining function $r$, we may write, in a
collar neighborhood of the infinity,
$$
g = \rho^{-2}(dr^2 + g_r)
$$
where $\rho = \sinh r$,
$$
g_r = h_0 + \frac {r^{n+1}}{n+1} h + O(x^{n+2})
$$
and $h$ is a symmetric 2-tensor on $S^n$.
\end{deff}

A function $u$ which is $m$-times continuously differentiable on
$X$ is said to be in the {\em weighted H$\ddot{o}$lder space}
$\Lambda^s_{m, \alpha}(X)$ if $|| u ||^s_{m, \alpha} < \infty$ for
$s\in R, m \geq 0$ and $\alpha \in (0,1)$, where the norm $|| u
||^s_{m, \alpha} $ is defined as follows. First, in the special
case in which $X$ is a smoothly bounded open subset of $\R^{n+1}$,
we define
$$
||u||^s_{m, 0} = \sum^m_{l=0} \sum_{|\gamma|=l}
||d^{-s+l}\partial^\gamma u||_{L^\infty}
$$
and
$$
||u||^s_{m, \alpha} =  || u ||^s_{m, 0} + \sum_{|\gamma|= m}
\sup_{x, y} \left[ \min (d^{-s + m + \alpha}_x, d^{-s + m + \alpha}_y)
\frac{ | \partial^\gamma u(x) - \partial^\gamma u(y) | }{|x-y|^\alpha}
\right] ,
$$
where $d_x$ is the Euclidean distance from $x$ to $\partial X$. In
the more general cases of a manifold with boundary, the same norms
are defined using a covering by coordinate charts and a
subordinate partition of unity in the usual way. We recommend
\cite{G_L} and \cite{Lee} for succinct discussions of properties
of the spaces $\Lambda^s_{m, \alpha}(X)$.

\section{Analytic and Geometric Preliminaries}

We first recall the following lemma from, for instance, \cite{G_L}
\cite{Lee}.

\begin{lma} Let $(X, g)$ be a weakly asymptotically hyperbolic manifold of
order $C^{3, \alpha}$. Then any representative $\hat{g}$ in the
conformal infinity of $g$ determines a unique defining function $s
\in C^{2, \alpha}(\bX)$ such that $s^2 g |_{TM} = \hat{g}$, $s^2g$
has a $C^{2, \alpha}$ extension to $\bar{X}$ and $|ds|^2_{s^2 g}
\equiv 1$ on a neighborhood $U$ of $M$ in $\bX$. Hence, $s$
determines an identification of $U$ with $M \times [0, \epsilon)$
such that
\begin{equation}
\label{snormal} g = \frac{1}{s^2} ( ds^2 + g_s )
\end{equation}
for a $1$-parameter family $\{g_s\}$ of metrics on $M$ with $g_0 = \hat{g}$.
\end{lma}

By a change of variable
\begin{equation}
s= \frac{\cosh(r) - 1}{\sinh(r)},
\end{equation}
we may rewrite (\ref{snormal}) as
\begin{equation}
\label{rhonormal}g = \rho^{-2} ( dr^2 + g_r ),
\end{equation}
where $\rho = \sinh(r)$. One may compare (\ref{rhonormal}) with
the fact that
$$
g_{b}=\frac{1}{\sinh^2(r)} \{ dr^2 + h_0 \}
$$
gives the standard hyperbolic metric on $S^n \times \R^{+}$ where
$h_0$ is the standard metric on $S^n$. The fact that $s$ is $C^{2,
\alpha}$ guarantees that the family of metrics $\{g_r\}$ is at
least $C^1$ with respect to $r$. In the special case in which
$(X^{n+1}, g)$ is Einstein and conformally compact of sufficiently
high order, Andersson and Dahl \cite{A_D} showed that the family
of metrics $\{g_r\}$ in (\ref{rhonormal}) have the properties
\begin{equation}
g_r = h_0 + \rho^n h, \ Tr_{h_0}h = O(\rho^n), \ \rho=\sinh(r).
\end{equation}
Thus, the decay assumption (\ref{decayassum}) is automatically
satisfied by any conformally compact Einstein manifold with the
round sphere as its conformal infinity. Next we recall an analytic
result of the operator $-\Lap_g + (n+1)$ between suitable weighted
H$\ddot{o}$lder spaces (see Proposition $3.3$ in \cite{Lee}).

\begin{lma} \label{isomorphism} Let $(X^{n+1}, g)$ be weakly asymptotically hyperbolic
of order $C^{m, \alpha}$. Let $0 < \beta < 1$ and $ k + 1 + \beta \leq m + \alpha$.
Then
$$
-\Lap + (n+1) : \Lambda^s_{k+2, \beta} \rightarrow \Lambda^s_{k, \beta}
$$
is an isomorphism whenever $ -1 < s < n+1$.
\end{lma}

In the final step of the proof in \cite{Qing}, the positive mass
theorem is used on the doubling of a partially compactified
manifold along its totally geodesic boundary. Here we observe that
it would be much simpler if we appeal to the following quasi-local
mass type result proved in \cite{Miao} \cite{S_T}(see also
\cite{Qing1}).

\begin{prop} \label{qlmass}
Let $\bar{\Omega}^{n+1}$ be a smooth compact manifold with boundary $ \partial \Omega$.
Let $g$ be a metric on $\bar{\Omega}$ which is smooth in the interior $\Omega$
and $C^2$ up to $\partial \Omega$. If $g$ has nonnegative scalar curvature
in $\Omega$, $(\partial \Omega, g|_{T \partial \Omega})$ is isometric to
$(S^n, h_0)$ and the mean curvature of $\partial \Omega$ with respect to
the outward pointing unit normal identically equals the constant $n$, then
$g$ has vanishing scalar curvature in $\Omega$ provided the dimension
satisfies $2 \leq n \leq 6$.
\end{prop}

\begin{rmk} It is desirable to further conclude that $g$ is actually flat on
$\Omega$ which is indeed the case when $n=2$ \cite{Miao}. However,
no proof in higher dimension is known so far except the case
when $(\Omega, g)$ is assumed to be spin \cite{S_T}.
\end{rmk}

We conclude this section by recalling a nice functional characterization of the
Hyperbolic space $\H^{n+1}$ proved in \cite{Qing1}.

\begin{lma} \label{hlma}
Let $(X^{n+1}, g)$ be a complete Riemannian manifold.
Assume that there exits a positive smooth function $u$ on $X$ such
that
$$
Hess_g (u) = u g.
$$
Then $(X^{n+1}, g)$ is isometric to $(\H^{n+1}, g_{\H})$.
\end{lma}

\section{Proof of the main Theorem}

Let $(X^{n+1}, g)$ satisfy the assumptions in Theorem
\ref{Rigidity} and let $U$ be a neighborhood of $M$ in $\bX$ where
(\ref{gexpan}) holds. We introduce a background hyperbolic metric
\begin{equation}
\label{gb}g_b = \frac{1}{\rho^2}\{ dr^2 + h_0 \}
\end{equation}
on $U$. Clearly $(U, g_b)$ can be identified with the complement of some compact set
in the Hyperbolic space $(\H^{n+1}, g_{\H})$ realized as the hypersurface
$$
\{ (x_1, \ldots, x_{n+1}, t) \ | \ |x|^2 - t^2 = -1, t>0 \} \subset
\R^{n+1, 1}
$$
by letting $\sinh r = \rho = \frac{1}{|x|}$. The restriction of
$t$ to $\H^{n+1}$ is an eigenfunction of $(\H^{n+1}, g_{\H})$,
i.e.
$$
\Delta_{\H} t = (n+1) t.
$$
Moreover, as observed in \cite{Qing}, by a change of variables,
$$
t = \frac {1+|y|^2}{1-|y|^2}
$$
and $(\H^{n+1}, g_{\H}) = (B^{n+1}, (\frac 2{1-|y|^2})^2|dy|^2)$,
which tells us $(t+1)^{-2} g_{\H}$ compactifies $\H^{n+1}$ to be
the standard Euclidean ball $\bar{B}^{n+1} \subset \R^{n+1}$ with
totally umbilical boundary $S^n$. This leads us to transplant $t$
to the domain $U$ and then look for a positive eigenfunction $u$
on $(X^{n+1}, g)$ which behaves like $t$ near $M$. To simplify
notations, we use $v \in O_k(\rho^s)$ to standard for $v \in
\Lambda^s_{k, \beta}(X)$ for $k \geq 0$ and a fixed $0 < \beta <
1$.

\begin{lma} \label{ugh} There exists a smooth function $u>0$ on
$(X, g)$ such that
\begin{equation}
-\Lap_g u + (n+1) u = 0
\end{equation}
and
\begin{equation} u = t + O_2(\rho^{\tilde{s}})
\end{equation}
for some $1 < \tilde{s} < n+1$.
\end{lma}
\pf The fact $r \in C^2(U)$ and $t= \sqrt{1+ \frac{1}{\rho^2}}$ implies
$t \in C^2(U)$. We calculate
\begin{eqnarray}
-\Lap_g t & = & - \frac{\rho^{n+1}}{\sqrt{det g_r}}
\partial_r (\rho^{1-n} \sqrt{det g_r} \partial_r t) \nonumber \\
& = & - (n+1)t + \frac{1}{2} Tr_{g_r}{g^\prime_r},
\end{eqnarray}
where ``$\prime$'' denotes differentiation with respect to $r$.
By the decay assumption (\ref{decayassum}), we may choose $1 < \tilde{s}
< \min(s, n+1)$ such that $Tr_{g_r}{g^\prime_r} = O_0(\rho^{\tilde{s}} ) .$
Hence, by Lemma \ref{isomorphism}, we know there exists a function
$w = O_2(\rho^{\tilde{s}})$ such that
\begin{equation} \label{defofu}
-\Lap_g w + (n+1) w = -\Lap_g t + (n+1)t .
\end{equation}
Let $u = t - w$. Then $u>0$ by the maximum principle and the
smoothness of $u$ follows directly from the local elliptic
regularity theory. \stop

We refer readers to \cite{Lee}, \cite{Qing} and \cite{Q_P_A}  for
more results on eigenfunctions for asymptotically hyperbolic
manifolds. In our next lemma we set the stage to apply the work
from \cite{Miao} by using the eigenfunctions to compacitify the
weakly asymptotically hyperbolic manifolds.

\begin{lma} The metric $g_u = \frac{1}{(1+u)^2} g$ extends to a $C^2$ metric
on $\bX$ such that $g_u$ has nonnegative scalar curvature in $X$,
$(M, g_u | _{TM})$ is isometric to $(S^n, h_0)$ and the mean curvature
of $M$ in $(\bX, g_u)$ identically equals the constant $n$.
\end{lma}
\pf First we calculate the scalar curvature of $g_u$,
\begin{eqnarray}
R(g_u) & = & \frac{4n}{1-n} (u+1)^{\frac{n+3}{2}}
[\Lap_g - \frac{n-1}{4n} R(g) ](u+1)^{\frac{1-n}{2}} \nonumber \\
& = & -n(n+1)|du|^2_g + 2n(n+1)u(u+1) + R(g)(u+1)^2  \label{defR} .
\end{eqnarray}
Since $Ric(g) \geq -n g$, we have $R(g) \geq -n(n+1)$ so
(\ref{defR}) implies
\begin{eqnarray}
R(g_u) & \geq &
-n(n+1)|du|^2_g + 2n(n+1)u(u+1) - n(n+1) (u+1)^2  \nonumber \\
& =
& n(n+1)( u^2 - |du|^2_g - 1 ) \label{Rgu} .
\end{eqnarray}
As in \cite{Qing}, we then appeal to the Bochner formula for
eigenfunctions, which is observed in \cite{Lee}.
\begin{eqnarray} \label{defB}
\Lap_g(|du|^2_g - u^2)
& = & 2n|du|^2_g + 2 Ric(\nabla_g u, \nabla_g u)
+ 2 |Hess_g
u|^2_g - 2(n+1)u^2  \nonumber \\
& \geq &  2 |Hess_g u|^2_g - 2(n+1)u^2 ,
\end{eqnarray}
where the last step holds again since $Ric(g) \geq -n g$.
Therefore, we have
\begin{equation}
\label{Bochner} - \Lap_g ( u^2 - |du|^2_g - 1 ) \geq 2 |Hess_g u -
u g|^2_g .
\end{equation}
Hence, in order to prove the scalar curvature $R(g_u) \geq 0$, we
only need to apply a maximum principle to $u^2 - |du|^2 - 1$ and
verify that it goes to zero towards the boundary. A
straightforward calculation reveals that
\begin{eqnarray}
u^2 - |du|^2_g & = & t^2 - 2tw + w^2 - \rho^2
(\partial_r t)^2 -
\rho^{2} 2\partial_r t \partial_r w
\nonumber
\\
& & - \rho^{2} (\partial_r w)^2 -
g^{\delta
\lambda}\partial_\delta w \partial_\lambda w \nonumber \\
& = & 1
+ O_2(\rho^{\tilde{s}-1}) + O_2(\rho^{2\tilde{s}})
+
O_1(\rho^{\tilde{s}-1}) \nonumber \\
& & +  O_1(\rho^{2\tilde{s}})
+ O_1(\rho^{2\tilde{s}}) ,
\end{eqnarray}
where we have used the fact $t^2 - \rho^2 (\partial_r t)^2 = 1$.
It follows from $\tilde{s}>1$ that
\begin{equation}
u^2 - |du|^2_g - 1  \rightarrow 0, \ \ as\
\rho \rightarrow 0 .
\end{equation}
Thus, we have
\begin{equation} \label{middle}
u^2 - |du|^2_g - 1 \geq 0 \ \
on\ X,
\end{equation}
which implies $R(g_u) \geq 0$ on $X$ by (\ref{Rgu}).

Next we consider the expansion of $g_u$ near $M$,
\begin{equation}
g_u = \frac{1}{[(u+1)\rho]^2}\{ dr^2 +
g_r\},
\end{equation}
where
\begin{equation}
(u+1)\rho = \cosh r +
\sinh r - w \sinh r .
\end{equation}
Since $w = O_2(\rho^{\tilde{s}}), \tilde{s} > 1$ and $\rho^2 g =
dr^2 + g_r$ is $C^2$ on $\bX$, we see that $g_u$ readily extends
to a $C^2$ metric on $\bX$. Furthermore, we have the boundary
values
\begin{equation}
(u+1)\rho |_{r=0} = 1 \ \ and\ \ \frac{d}{dr} [(u+1)\rho] |_{r=0}
= 1 ,
\end{equation}
which, combined with the facts $g_0 = h_0$ and $Tr_{g_0}g^\prime_0
= 0$, show that $(M, g_u | TM)$ is isometric to $(S^n, h_0)$ and
$M$ has constant mean curvature $n$ in $(\bX, g_u)$.

Now it follows from Proposition \ref{qlmass} that $R(g_u) \equiv
0$ on $X$. (\ref{Rgu}), (\ref{middle}) and (\ref{Bochner}) then
imply that
\begin{equation}
|Hess_g u - u g | =0.
\end{equation}
Therefore $(X^{n+1}, g)$ is the hyperbolic space $\H^{n+1}$ by
Lemma \ref{hlma}. \stop

\begin{rmk}
We may reformulate the decay assumption (\ref{decayassum}) in
terms of the metric expansion (\ref{snormal}). By substituting $s
= \frac{\cosh (r) - 1}{\sinh (r)}$ back we see that
(\ref{decayassum}) is equivalent to
\begin{equation}
(1 -
\frac{s^2}{4}) Tr_{g_s} (\frac{d}{ds} g_s ) + ns \in
\Lambda^{\delta}_{0, \beta}(X)
\end{equation}
for some $\delta > 1$.
\end{rmk}

To conclude we would like to make some remarks. The main theorem
in this paper improves the rigidity theorem in \cite{Qing}. Here
we no longer assume that the conformally compact manifolds are
Einstein and we assume much weaker asymptotics at the infinity. In
other words, we only assume the Einstein equations are satisfied
at the infinity to a very low order, which is often true since the
energy-momentum tensor usually vanishes to certain order for
isolated systems. Also, we believe it is interesting to compare
our result to the scalar curvature rigidity for asymptotically
hyperbolic manifolds in \cite{A_D}, \cite{C_H}, \cite{Min-Oo},
\cite{Zhang}. For example, in \cite{Wang}, Wang defines a
conformally compact manifold $(X^{n+1}, g)$ to be {\em
asymptotically hyperbolic} if it satisfies:
\begin{enumerate}

\item $(X^{n+1}, g)$ is weakly asymptotically hyperbolic with the conformal
infinity being the
standard sphere $(S^n, h_0)$.

\item Let $r$ be the special
defining function so that we can write
\begin{equation}
g =
\frac{1}{\sinh^2 (r)} \{ dr^2 + g_r \}
\end{equation}
in a neighborhood of $\partial X$. Then
\begin{equation}
\label{wdecay}
g_r = h_0 + \frac{r^{n+1}}{n+1}h + O(r^{n+2}),
\end{equation}
where $h$ is a symmetric $2$-tensor on $S^n$. Moreover the
asymptotic expansion can be differentiated twice.
\end{enumerate}

Working with this definition, see also \cite{A_D}, \cite{C_H} and
\cite{Min-Oo}, Wang was able to prove that if $(X^{n+1}, g)$ is
asymptotically hyperbolic, $(X^{n+1}, g)$ is spin and the scalar
curvature $R \geq -n(n+1)$, then
\begin{equation*}
\int_{S^n} ( Tr_{h_0} h ) d\mu_{h_0} \geq \left| \int_{S^n} (
Tr_{h_0} h ) x d\mu_{h_0} \right|.
\end{equation*}
Moreover equality holds if and only if $(X, g)$ is isometric to
the hyperbolic space $\H^{n+1}$. Since our decay assumption in
Theorem \ref{Rigidity} is much weaker than (\ref{wdecay}), we
immediately have the following corollary,

\begin{corol}
Let $(X^{n+1}, g)$ be an asymptotically hyperbolic manifold in the
sense of \cite{Wang}. If $ 2 \leq n \leq 6$ and $Ric \geq -n g$,
then $(X^{n+1}, g)$ is isometric to the hyperbolic space
$\H^{n+1}$.
\end{corol}

\bibliographystyle{plain}
\bibliography{rigidity}

\end{document}